\documentclass[11pt]{article}
\usepackage[margin=1in]{geometry}
\usepackage{amsmath,amsfonts,amssymb,amsthm}
\usepackage{graphicx}
\usepackage{enumerate}
\usepackage{bbm}
\usepackage{verbatim}
\usepackage{hyperref,color}
\usepackage[capitalize,nameinlink]{cleveref}
\usepackage[dvipsnames]{xcolor}
\hypersetup{
	colorlinks=true,
	pdfpagemode=UseNone,
    citecolor=OliveGreen,
    linkcolor=NavyBlue,
    urlcolor=Magenta,
	pdfstartview=FitW
}
\usepackage{appendix}
\crefname{appsec}{Appendix}{Appendices}
\usepackage{tikz}
\usepackage[colorinlistoftodos]{todonotes}
\usepackage{fourier}

\theoremstyle{plain}
\newtheorem{thm}{Theorem}[section]
\newtheorem{prop}[thm]{Proposition}

\newtheorem{lemma}[thm]{Lemma}

\newtheorem{conj}[thm]{Conjecture}

\newtheorem{qn}[thm]{Question}

\theoremstyle{definition}
\newtheorem{defn}[thm]{Definition}

\newtheorem{obs}[thm]{Observation}

\newtheorem*{ass*}{Assumption}

\theoremstyle{remark}

\newtheorem{rem}[thm]{Remark}

\crefname{lem}{Lemma}{Lemmas}
\crefname{thm}{Theorem}{Theorems}
\crefname{defn}{Definition}{Definitions}
\crefname{fact}{Fact}{Facts}
\crefname{clm}{Claim}{Claims}
\crefname{prop}{Proposition}{Propositions}

\newcommand{\rand}{\overset{\text{R}}\leftarrow}

\newcommand{\mc}{\mathcal}

\newcommand{\FF}{\mathbb{F}}

\newcommand{\R}{\mathbb{R}}
\newcommand{\RR}{\mathbb{R}}

\newcommand{\EE}{\mathbb{E}}
\newcommand{\ZZ}{\mathbb{Z}}

\newcommand{\mcF}{\mathcal{F}}

\newcommand{\mcS}{\mathcal{S}}

\newcommand{\ind}{\mathbf{1}}

\renewcommand{\P}{\mathbb{P}}

\begin{document}
	
\title{On sum-intersecting families of positive integers}
\author{Aaron Berger and Nitya Mani (MIT)}
\date{\today}

\maketitle

\begin{abstract}
We study the following natural arithmetic question regarding intersecting families: how large can a family of subsets of integers from $\{1, \ldots n\}$ be such that, for every pair of subsets in the family, the intersection contains a \textit{sum} $x + y = z$? We conjecture that any such \textit{sum-intersecting} family must have size at most $\frac14 \cdot 2^{n}$ (which would be tight if correct). Towards this conjecture, we show that every sum-intersecting family has at most $0.32 \cdot 2^n$ subsets.
\end{abstract}

\section{Introduction}
The study of intersection problems for finite sets and graphs begins with the following classical problem: determine the maximum size for a family $\mcF$ of subsets of $[n] = \{1, 2, \ldots , n\}$ such that for any two sets $A, B \in \mcF$, the intersection $A \cap B$ is nonempty. In other words, what is the maximum size of an \textit{intersecting family} of subsets of $[n]$? By taking $\mcF$ to be the collection of subsets of $[n]$ containing a fixed element, we can find an intersecting family of size $\frac12 \cdot 2^n$, and this bound is tight, as for any $A \subseteq [n]$, at most one of $A, A^c$ can be in any intersecting family.

More broadly, unstructured intersecting families of subsets of $[n]$ are generally well understood. For example, the Erd\H{o}s-Ko-Rado theorem characterizes maximal intersecting families of subsets with a fixed size: the theorem states that if $\mcF$ is an intersecting family of subsets of size $r \le n/2$, then $|\mcF| \le {n - 1 \choose r-1}$, a bound achieved by the collection of size-$r$ subsets of $[n]$ that contain a fixed element. The next question one might ask is for families in which the intersection of any two subsets contains some \textit{arithmetic structure}. For example, in 1986, Chung, Graham, Frankl, and Shearer~\cite{CGFS86} conjectured that the maximum size of a \textit{nontrivial $3$-AP intersecting} family of subsets of $[n]$ is $\frac18 2^n$ (i.e. a family $\mcF \subset 2^{[n]}$ where for every $A, B \in \mcF$, $A \cap B$ contains $x, y, z$ not all equal with $2y = x + z$); this conjecture remains wide open. In the same paper, they characterize maximal families where subsets must intersect on $k$ consecutive integers; the $k = 2$ case will serve as a motivating example for \cref{c:upper} to follow. Intersection problems concerning translates of arbitrary fixed patterns have received further attention; see for instance, the surveys~\cite{DEZ83,BOR12,ELL22}. 

We should remark that one may impose other types of structure on the ground set as well. Examples include intersecting families of combinatorial objects like permutations~\cite{EFP11}, and $H$-intersecting families of graphs on $n$ labelled vertices (where $H$ is a fixed graph, such as a triangle~\cite{EFF12} or larger clique~\cite{BZ21}).

Many known methods for proving intersecting family results are broadly applicable to structures which are \textit{avoided} by a randomized family of dense sets; for example, the random complete bipartite graph plays a central role in the triangle-intersecting family problem. This is one factor behind the difficulty of the 3-AP problem, as by Roth's theorem no dense set avoids 3-APs. By contrast, there is a classical construction of dense sum-free subsets of sets of integers, and tight bounds on the maximal sum-free subset problem are a topic of much interest \cite{BOU97,EGM14,SHA22}. We therefore believe that sum-intersecting families are a natural topic for study and represent an obvious gap in the literature; in this article we begin to fill this gap.

\begin{defn}
A collection $\mcF$ of subsets of $\{1, 2, \ldots, n\}$ is \textit{sum-intersecting} if for every $A, B \in \mcF$, $A \cap B$ contains a \textit{sum}, i.e. there exist $x, y, z \in A \cap B$ (not necessarily distinct) with $x + y = z$. 
\end{defn}

Note that for every sum-intersecting family $\mcF$, we have $|\mcF| \le \frac12 2^n$, since for any $A \subset [n]$, at most one of $A, A^c \in \mcF$. Further, note that the collection $\{S \subseteq [n] \mid \{1, 2\} \subseteq S\}$ is sum-intersecting, giving a sum-intersecting family with size $\frac14 2^n$. In this work we begin to close bring these bounds closer together.

The $\{x,x+1\}$-intersecting result of \cite{CGFS86} has implications for a finite field analogue of the problem we study here. Specifically, if $2$ is a generator of the cyclic group $\FF_p^{\times}$ for odd prime $p$, then~\cite{CGFS86} shows that the largest size of a $\{2x = y\}$-intersecting family in $\FF_p^{\times}$ is $\frac14 2^{p-1}.$ This shows that the construction $S \supset \{1,2\}$ is optimal for $\{2x = y\}$-intersecting families over $\FF_p$. This result motivates us to conjecture the same is true for all \textit{sum}-intersecting families over $[n]$.

\begin{conj}\label{c:upper}
If $\mcF$ is a sum-intersecting family of subsets of $\{1, 2, \ldots, n\}$, then $|\mcF| \le \frac14 2^{n}.$
\end{conj}

We tackle this question with a Fourier analytic reduction akin to the method used by Ellis, Filmus, and Friedgut~\cite{EFF12} to study triangle-intersecting families of graphs. Specifically, we use a linear programming bound (\cref{lem:lp_bound}) in conjunction with a study of distribution of random sum-free subsets of integers. Our methods give a recipe for the reduction of~\cref{c:upper} to a finite (although not clearly tractable) computation. Thus, while we do not completely resolve~\cref{c:upper}, we are able to substantially improve upon a simple entropy upper bound. As a specific illustration of our methods, we show in this work the following upper bound on the maximum size of a sum-intersecting family:

\begin{thm}\label{t:main}
If $\mcF$ is a sum-intersecting family of subsets of $\{1, 2, \ldots, n\}$, then $|\mcF| \le \frac{8}{25} \cdot 2^{n} = 0.32 \cdot 2^n.$
\end{thm}
We emphasize that we have optimized our presentation for simplicity and clarity here; at the expense of more complicated arguments along the same vein and increasingly large computations, the above bound can be improved further. However, we believe some new ideas will be necessary to reduce the full verification of \cref{c:upper}, if it is true, from a finite computation to a \textit{tractable} finite computation.
\newline \newline 
\textbf{Acknowledgements:} Thank you to Ashwin Sah for helpful discussions regarding~\cite{BOU97}. Thank you to Yufei Zhao for helpful comments on a draft. NM was supported by a Hertz Graduate Fellowship. NM and AB were supported by the NSF GRFP DGE-1745302.
\newline \newline 
\textbf{Notation:} We let $[n] = \{1, 2, \ldots, n\}$ and use $\alpha \rand \Omega$ to denote a sample from the uniform distribution on $\Omega.$ We let $\{z\} = z - \lfloor z \rfloor$ be the \textit{fractional part of} $z$. For $\alpha \in [0, 1]$ and $S \subseteq [n],$ we let $S_{\alpha} = \{x \in S \mid \{\alpha x\} \in [1/3, 2/3)\}.$

\section{Upper bounds}
Any sum-intersecting family has size at most $\frac 12 2^{n}$, as it is intersecting. One can obtain a first improvement to this trivial bound by applying Shearer's Lemma, a powerful tool from entropy first introduced in \cite{CGFS86} to tackle the triangle-intersecting family problem. An analogous argument to that of \cite{CGFS86} yields a bound of $2^{n-3/2} < 0.354 \cdot 2^n$. We encourage the reader to keep this bound in mind as a point of comparison for what follows, but we will omit further details as this bound will be superseded by the arguments we present next, and the entropy method will not be relevant there.

To improve on the entropy bound, we leverage the framework of~\cite{EFF12} to give a Delsarte-type linear programming bound on the maximum size of an intersecting family. 
Practically speaking, this framework allows us to give an upper bound on the size of an intersecting family by verifying simple linear inequalities regarding the probabilities of various outcomes when one takes a random sum-free subset of an arbitrary set of integers.

\begin{defn}
For a set $S \subset [n],$ we let $S_{\alpha}$ be the random set $\{x \in S \mid \{\alpha x \} \in [1/3 , 2/3)\}$ for $\alpha \rand [0, 1].$
\end{defn}
Observe that $S_\alpha$ is always sum-free.
\begin{lemma}[LP bound]\label{lem:lp_bound}
    Fix $c_0, c_1 \in \R$, and let $\mu : 2^{[n]} \to \R$ be given by $\mu(S) = (-1)^{|S|}(c_0 \P[|S_\alpha| = 0] + c_1 \P[|S_\alpha| = 1])$. Then if $\mu(S) \ge -1$ for all $S$, any sum-intersecting family $\mc F$ of subsets of $[n]$ has size at most $\frac 1{1+c_0}2^n$.
\end{lemma}
Thus, our goal will be to choose $c_0, c_1$ to maximize $c_0$ subject to the infinitely many (linear in $c_0, c_1$) constraints $\mu(S) \ge -1$ for all $S \subset [n]$; this motivates our characterization of the above as a special case of a \textit{linear programming bound}. We derive~\cref{lem:lp_bound} at the end of~\cref{s:lp-red} as a consequence of a more general result. 

\subsection{LP reduction}\label{s:lp-red}
Below we employ the following standard conventions for Fourier analysis over $\FF_2^n$: given $f: \FF_2^N \to \RR$, its \textit{Fourier transform} is given by 
$$\widehat{f}(\lambda) = \EE_{x \rand \FF_2^n}\left[(-1)^{\langle \lambda, x \rangle} f(x)\right], \quad \lambda \in \FF_2^n,$$
and we have \textit{Fourier inversion formula}
$$f(x) = \sum_{\lambda \in \FF_2^n} (-1)^{\langle \lambda, x \rangle} \widehat{f}(\lambda), \quad x \in \FF_2^n.$$
The \textit{convolution} of two functions $f, g: \FF_2^n \to \RR$ is given by $f * g: \FF_2^n \to \RR$ with $f * g(z) = \EE_{y \rand \FF_2^n}[f(y)g(z -y)]$. We have that $\widehat{f * g}(\lambda) = \widehat{f}(\lambda) \widehat{g}(\lambda).$ We will derive our upper bound from the following LP bound, which is essentially a simpler version of Proposition 2.1 of~\cite{BZ21}.
\begin{prop}\label{p:fourier}
Let $\mcF$ be a family of subsets of $[n]$ with associated indicator function $\ind_{\mcF}: \FF_2^n \to \{0, 1\}$. Suppose we are given $\nu: \FF_2^n \to \RR$ such that $\hat \nu \ge -1$ and $(\ind_{\mcF} * \ind_{\mcF} * \nu)(0)  = 0$. Then  
$$\frac{|\mcF|}{2^n} \le \frac{1}{1+\hat \nu(0)}.
$$
\end{prop}
\begin{proof} We give a streamlined proof here; this method appears in more detail in \cite{BZ21,EFF12}. Write
\begin{equation*}
    0 = f*f*\nu(0) = \sum_\lambda \widehat{f * f * \nu}(\lambda) = \sum_\lambda \hat f(\lambda)^2\hat\nu(\lambda) = \hat f(0)^2\hat\nu(0) + \sum_{\lambda \neq 0} \hat f(\lambda)^2 \hat \nu(\lambda) \ge \hat f(0)^2\hat\nu(0) - (\hat f(0) - \hat f(0)^2).
\end{equation*}
    Solving for $\hat f(0) = |\mcF|2^{-n}$ yields the proposition.
\end{proof}
In order to apply~\cref{p:fourier}, we will let $\mcF$ be an arbitrary sum-intersecting family and we will need to choose a suitable $\nu$. 
Constructing a large sum-intersecting family $\mc F$ can be viewed as solving a particular linear program; constructing a suitable $\nu$ to obtain a bound on a maximal $\mc F$ is simply solving the dual LP. Methods to construct suitable $\nu$ are now well-understood. The following result essentially appears in  \cite{BZ21} and \cite{EFF12}.
\begin{lemma}[Lemma 2.3~\cite{EFF12}, Lemmas 2.2 and 2.3~\cite{BZ21}] \label{lem:struct-nu}
Consider an arbitrary distribution $\mcS$ on sum-free subsets of $[n]$ and functions $\{f_S\}$ indexed by subsets of $[n]$, and define function $\nu: \FF_2^n \to \RR$ by its Fourier transform
\begin{equation}\label{eqn:def_nu}
    \widehat{\nu}(\lambda) = (-1)^{\|\lambda\|_1} \EE_{S \sim \mcS} f_S(S \cap \lambda).
\end{equation}
Then, for any sum-intersecting family $\mcF$,  $(\ind_\mcF * \ind_\mcF * \nu)(0) = 0$.
\end{lemma}
We present a streamlined proof here; we refer the reader to the given references for further details. 
\begin{proof}
    First, we verify the convolution identity when $\nu = \ind_{T}$ is a point mass at $T \subseteq [n]$ such that $S := [n] \setminus T$ is sum-free. In fact, a stronger result holds: $\ind_\mcF(S_1)\ind_\mcF(S_2)\ind_T(S_3) = 0$ whenever $S_1+S_2+S_3 = 0$; otherwise, we would need $S_3 = T$, and so $S_2 = S_1+T$. We know $S_1,S_2 \in \mcF$, but $S_1 \cap S_2 \subseteq [n] \setminus T$ is sum-free, which contradicts $\mcF$ being sum-intersecting.

    The Fourier transform of this $\nu$ can be written in the form of (\ref{eqn:def_nu}): one can verify that taking $\mcS$ to be a point mass at $S$ and $f_S(X) = (-1)^{|X|}$ does indeed yield $\hat \nu$. The set of all such $\hat \nu$ for $S' \subseteq S$ forms a basis for the space of all $f_S$, hence (\ref{eqn:def_nu}) holds for all $f_S$ by linearity, and then for all $\mcS$ by applying linearity again.
\end{proof}

\cref{lem:lp_bound} arises by specializing this general result to a simple distribution $\mcS$ and choice of $f_S$.

\begin{proof}[Proof of~\cref{lem:lp_bound}]
Given $c_0, c_1 \in \RR$, we construct $\nu: \FF_2^n \to \RR$ as function arising from choice of Fourier transform in~\cref{lem:struct-nu} as follows. In the notation of~\cref{lem:struct-nu}, we define
$$\mcS \sim [n]_\alpha : \alpha \rand [0,1] \quad \text{ and }\quad f_S(T) := c_0 \ind(T = \emptyset) + c_1 \ind(|T| = 1).$$

This yields $\hat \nu  = \mu$. By construction $\hat \nu(\emptyset) = c_0$ and $\hat \nu \ge -1$ always. By \cref{lem:struct-nu} $\nu$ satisfies $ (\ind_{\mcF} * \ind_{\mcF}* \nu) (0) = 0$ for any sum-intersecting family $\mc F$. Thus, $\nu$ satisfies the hypotheses of \cref{p:fourier}, and we conclude $|\mcF| \le \frac{1}{1+c_0}2^{-n}$, completing the proof.
\end{proof}

Thus, in the remainder of the article, our goal will be to find the maximal choice of $c_0 \in \RR$ such that~\cref{lem:lp_bound} holds for some $c_1 \in \RR$.

\subsection{A $\frac13$ upper bound}
Recall that for a set $S \subset [n],$ we let $S_{\alpha}$ be the random set $\{x \in S \mid \{\alpha x \} \in [1/3 , 2/3)\}$ for $\alpha \rand [0, 1].$
We require the following result of Shakan, which gives an exact description of the distribution of $|S_\alpha|$ when $|S| = 2$.
\begin{thm}[Proposition 1~\cite{SHA22}]\label{thm:2_point_formula}
    Let $S = \{x, y\} \subset [n]$. Define
    \begin{equation}
    E(x, y) := \frac{\chi(xy) \gcd(x, y)^2}{xy},
    \end{equation}
    where $\chi(k)$ is the multiplicative character mod $3$ given by $\chi(0) = 0, \chi(1) = 1, \chi(2) = -1.$ Then
    $$\P[|S_\alpha| = 2] = \frac{1}{9} + \frac29 E(x, y).$$
\end{thm}

We call $E$ the \textit{error function} because it is 0 precisely when the events $x \in S_\alpha$ and $y \in S_\alpha$ are independent, and otherwise it is directly proportional to the discrepancy between this idealized independent case and the truth. We begin by making a few simple observations about $E(a, b)$ that follow directly from the fact that for $a,b\in [n]$, $E(a,b) = E\left(\frac{a}{gcd(a,b)}, \frac{b}{gcd(a,b)}\right)$.
\begin{obs}\label{o:error}
The following properties are enjoyed by $E(a, b)$ for $a\neq b \in \ZZ_{> 0}$:
\begin{itemize}
    \item $E(a, b) \le \frac14,$ which is achieved when $a/b \in \{1/4, 4\}.$
    \item If $E(a, b) > 0,$ then $E(a, b) = \frac{1}{3k+1}$ for some $k \in \ZZ_{>0}.$
\end{itemize}
\end{obs}

\begin{defn}
For $S = \{a_1, \ldots, a_k\}$ and a subset $\{i_1, \ldots, i_j\} \subset [k]$ for $j \ge 1$, we let 
$$p_{a_{i_1}\cdots a_{i_j}} = \P_{\alpha}\left(S_{\alpha} \supset \{a_{i_1}, \ldots, a_{i_j}\}\right).$$
We let $p_{\emptyset} = \P(S_{\alpha} = \emptyset)$ and $p_{\ind} = \P(|S_{\alpha}| = 1).$
\end{defn}

\begin{prop}\label{p:13}
If $\mcF$ is a sum-intersecting family of subsets of $\{1, 2, \ldots, n\}$, then $|\mcF| \le \frac13 \cdot 2^{n}$.
\end{prop}
\begin{proof}
Apply \cref{lem:lp_bound} with $c_0 = 2$ and $c_1 = -1$. We need to show that for all sets $S$, $\mu(S) = (-1)^{|S|}(2 \P[|S_\alpha| = 0] - \P[|S_\alpha| = 1]) \ge -1$. We explicitly verify $\mu(\emptyset) = 2$ and $\mu(S) = -1$ whenever $|S| = 1$, so we turn to the case $|S| \ge 2$. 

As probabilities lie in $[0,1]$, we see that $\mu(S) \ge -\P(|S_{\alpha}| = 1) \ge -1$ when $|S|$ is even and $\mu(S) \ge -2p_\emptyset$ when $|S|$ is odd. Therefore is enough to show $p_{\emptyset} \le \frac 12$ whenever $|S|$ odd and $|S| \ge 2$; in fact this holds for all $|S| \ge 2$. Since $p_{\emptyset}$ cannot increase when an element is added to $S$ it suffices to verify $p_{\emptyset} \le \frac12$ when $S = \{x, y\}$, which holds by the principle of inclusion and exclusion:
$$p_{\emptyset} = 1 - p_x - p_y + p_{xy} = 1 - \frac 13 - \frac13 + \left(\frac19 + \frac29 E(x, y)\right) \le \frac13 + \frac19 + \frac29 \cdot \frac14 = \frac12,$$
where we use the fact that for any $x \in [n],$ $p_x = \frac13$, apply~\cref{thm:2_point_formula}, and in the final inequality we apply~\cref{o:error}. This gives the desired result.
\end{proof}

\subsection{Proof of~\cref{t:main}}\label{s:mainthmpf}

The proof of~\cref{p:13} reduced to the observation that $p_{\emptyset} \le \frac12$. In fact if $|S| \ge 3$, we can obtain a stronger upper bound that will yield a correspondingly stronger upper bound on the size of a maximum sum-intersecting family via~\cref{lem:lp_bound}. We employ the same notation as the previous subsection.

\begin{lemma}\label{l:1124}
If $S \subset [n],$ with $|S| \ge 3$, then $p_{\emptyset} \le \frac{11}{24}$. 
\end{lemma}
\begin{proof}
As before, we assume without loss of generality that $|S| = 3$ and write $S = \{x, y, z\}.$ By the principle of inclusion-exclusion,
$$p_{\emptyset} = 1 - p_x - p_y - p_z + p_{xy} + p_{yz} + p_{xz} - p_{xyz} \le p_{xy} + p_{yz} + p_{xz},$$
noting that $p_x = p_y = p_z = \frac13$ and $p_{xyz} \ge 0.$
We then apply~\cref{thm:2_point_formula} to find that
$$p_{\emptyset} \le \frac13 + \frac29(E(x, y) + E(y, z) + E(x, z)).$$
Thus, to show $p_{\emptyset} \le 
\frac{11}{24}$,
it suffices to show that $E(x, y) + E(y, z) + E(x, z) \le \frac{9}{16}.$ 

Suppose without loss of generality that $E(x, y) \ge E(x, z) \ge E(y, z)$. By \cref{o:error}, the value of $E(x,y)$, if it is positive, lies in $\{\frac 14, \frac 17, \ldots\}$ and depends only on $\frac xy$. Indeed, if $E(x,y) = E(x,z) = \frac 14$, then from the definition of $E$ we see that $\frac xy, \frac xz \in \{\frac 14, 4\}$. As $y \neq z$, this can only occur when $\{x,y,z\} \cong \{1,4,16\}$ up to permuting and rescaling, in which case we have $E(y,z) = \frac 1{16}$ and $E(x, y) + E(y, z) + E(x, z) = \frac{9}{16}$ as desired. 
Otherwise, we still have $E(x,y) \le \frac 14$, but now $E(y,z) \le E(x,z) \le \frac 17$. Then $E(x, y) + E(x, z) + E(y, z) \le \frac 14+\frac 17+\frac 17 = \frac{15}{28} < \frac 9{16}$ and we are done.
\end{proof}
\color{black}

\begin{rem}
The optimal bound for $|S| \ge 3$ would be to show $p_{\emptyset} \le \frac{7}{18}$, which is tight in the case $S = \{1,3,4\}.$ We omit a proof here as we do not require this stronger statement, but the argument begins similarly: the desired bound holds whenever $E(x,y) + E(x,z) + E(y,z) \le \frac 14$, which is true unless $x/y \in \{\frac 1{10}, \frac 25, \frac 17, \frac 14\}$, or their reciprocals. In these exceptional  cases, one would need to check the conclusion of the lemma directly, for which some extra work is required.
\end{rem}

We will be able to eke out a further improvement by obtaining an improvement over the trivial bounds $0 \le \P(|S_{\alpha}| = 1) \le 1$ when $|S| \ge 3$.
\begin{lemma}\label{l:bohr}
For $S \subset [n]$, we have that $p_{\emptyset} \ge \frac{2}{3^{|S|} - 2^{|S|} + 2}$.
\end{lemma}

We will prove~\cref{l:bohr} using the language of \textit{Bohr sets}, noting that for fixed $\alpha \in [0, 1]$, the subset $S \backslash S_{\alpha}$ are those elements $x \in S$ where $\alpha x$ has distance at most $1/3$ to the nearest integer. 

\begin{defn}
Let $G$ be a compact abelian group and for $z \in [0,1)$, let $\| z\| = \min(z, 1-z)$ be the distance from $z$ to the nearest integer.
Let $S  = \{a_1, \ldots, a_k\} \subseteq \hom(G, \RR/\ZZ)$ and choose $\rho \in [0, 1)$. The \textit{Bohr set} $B(S, \rho)$ is given by $$B(S, \rho) := \left\{x \in \RR/\ZZ : \|a_i x\| \le \rho \text{ for all } i \in [k]\right\}.$$
\end{defn}
\begin{rem}
In fact, in this work, we will only need to consider the case $G = \RR/\ZZ$ and Bohr sets $B(S, 1/3)$ for some $S \subseteq [n]$, where $a \in S$ corresponds to the automorphism $x \mapsto ax$. Note that $B(S, 1/3) = \{ \alpha : |S_{\alpha}| = 0\}$.
\end{rem}
The following is a standard lower bound on the size of Bohr sets.

\begin{lemma}\label{l:bohr-1}
Let $G$ be a compact abelian group, and $\rho \ge 1/t$ for some $t \in \ZZ_{> 0}$. Then
$|B(S, \rho)| \ge t^{-|S|}$.
\end{lemma}
\begin{proof}
Let $B = B(S, \rho)$ for $\rho \ge 1/t$ with $S = \{a_1, \ldots, a_k\}$. We will show that $t^{|S|}$ translates of $B$ cover $G$, from which the result follows. To see this, partition $[0, 1]$ into $t^{|S|}$ subsets of the form
$$Y_{b_1, \ldots, b_{k}} = \left\{x : \|a_i x\| \in \bigg[\frac{b_i}{t}, \frac{b_i + 1}{t}\bigg) \text{ for all } i \in [k]\right\},$$
where each $b_j \in \{0,1, \ldots, t-1\}$. We will show that for each such subset $Y = Y_{b_1, \ldots, b_{k}}$, we can cover $Y$ by a translate of $B$. Take $x, x' \in Y$ and notice that by the triangle inequality
$$\|a_i(x - x')\| = \|a_i x - a_i x'\| \in (-\rho, \rho),$$
and thus $x - x' \in B$ or that $x \in x' + B$ for all $x \in Y$. Thus $Y \subset x' + B$.
Therefore,
$$1 = \sum_{b_1, \ldots, b_k \in \{0, \ldots, t-1\}} |Y_{b_1, \ldots, b_k}| \le \sum_{b_1, \ldots, b_n \in \{0, \ldots, t-1\}} |B| = t^{k} |B|,$$
and thus as desired, $|B| \ge \rho^{|S|}$.
\end{proof}

If $G$ is an infinite group we can strengthen the above argument by slightly more than a factor of 2.
\begin{lemma}\label{l:strong-bohr}
Let $G$ be an infinite compact abelian group, and $\rho \ge 1/t$ for some $t \in \ZZ_{> 0}$. Then
$|B(S, \rho)| \ge \frac{2}{t^{k} - 2^k + 2}$.
\end{lemma}
\begin{proof}
For $B = B(S, \rho),$ let $B^+ = \{x : \{a_1 x\} \in [0, \rho]\} \cap B$ and let $B^- = B \backslash B^+$. 
Observe that $-B^+ = B^-$ up to a boundary of measure 0 (here we use the fact that $G$ is infinite), so $|B^+| = |B^-| = |B|/2$. Consider the argument of~\cref{l:bohr-1}, but let $x' \in Y$ be chosen as $\text{argmin}_{x \in Y} \|a_1 x\|$. Then, $x - x' \in B^+$ for all $x \in Y$, so in fact we can cover $Y$ by a translate of $B^+$. Further, note that $B^+, B^-$ collectively cover $2^k$ total subsets of the form $Y_{b_1, \ldots, b_k}$ (as defined in the proof of~\cref{l:bohr-1}). This gives a cover of $[0, 1]$ by at most $t^k - 2^k + 2$ translates/reflections of $B^+$ up to a measure $0$ boundary, showing that $|B| = 2|B^+| \ge \frac{2}{t^k - 2^k + 2}$. 
\end{proof}

\begin{proof}[Proof of~\cref{l:bohr}]
Note that $p_{\emptyset} = |B(S, 1/3)|$. The lower bound on $p_{\emptyset}$ immediately follows by applying~\cref{l:strong-bohr} with $\rho = 1/3$.
\end{proof}

\begin{proof}[Proof of~\cref{t:main}]
We apply \cref{lem:lp_bound} with $c_0 = 17/8$ and $c_1 = -5/4$. Then, it suffices to show that for all sets $S$, $$\mu(S) = (-1)^{|S|}\Big(c_0 \P(|S_{\alpha}| = 0) + c_1 \P(|S_{\alpha}| = 1)\Big) = (-1)^{|S|}\left( \frac{17}{8} p_{\emptyset} - \frac{5}{4}p_{\ind} \right) \ge -1.$$
We proceed by casework on $|S|$.
\begin{itemize}
    \item When $|S| = 1$, $p_{\emptyset} = 2/3, p_\ind = 1/3$; so $\mu(S) = - \frac{17}{8}\cdot \frac23 + \frac{5}{4} \cdot \frac13 = -1$ and thus the desired inequality holds.
    \item When $S = \{x, y\}$, observe $p_\ind \le p_x + p_y = 2/3$. Then $\mu(S) = \frac{17}{8} p_{\emptyset} - \frac{5}{4}p_{\ind} \ge - \frac{5}{4}\cdot\frac 23  = -\frac 56 \ge -1$.
    \item When $|S|$ is odd and at least $3$, we have $p_{\emptyset} \le 11/24$ by~\cref{l:1124}, and so
    $$\mu(S) = -\left(\frac{17}{8} p_{\emptyset} - \frac54 p_{\ind} \right) 
    \ge - \frac{17}{8} \cdot \frac{11}{24} = -\frac{187}{192} \ge -1.$$
    \item When $|S| = 4$, by~\cref{l:bohr}, $p_{\emptyset} \ge 2/67$. Further, observe that 
    $$\frac83 = \EE_{\alpha}[4 - |S_{\alpha}|] = \P(|S_{\alpha}| = 3) + 2 \P(|S_{\alpha}| = 2) + 3\P(|S_{\alpha}| = 1) + 4\P(|S_{\alpha}| = 0) \ge 3p_{\ind} + 4 p_{\emptyset} \ge 3p_{\ind} + \frac{8}{67}.$$
    Solving, we see that $p_{\ind} \le \frac{512}{603}$ when $|S| = 4$. Then, we can check that as desired,
    $$\mu(S) = \frac{17}{8} p_{\emptyset} - \frac54 p_{\ind} \ge \frac{17}{8}\frac{2}{67} - \frac54 \frac{512}{603} = -\frac{2407}{2412}\ge -1.$$
    \item When $|S|$ is even and at least $6$, then we observe 
    \begin{equation*}
        \frac{2|S|}{3} = \EE_{\alpha}[|S| - |S_{\alpha}|] \ge \sum_i (|S| - i)\P[|S_\alpha| = i] \ge (|S| - 1)p_\ind.
    \end{equation*}
    We conclude $P_\ind \le \frac{2}{3}\cdot\frac{|S|}{|S| - 1} \le \frac 23 \frac 65 = \frac 45$ when $|S| \ge 6$. Then the only negative contribution to $\mu$ is $c_0p_\ind \ge -1$.
\end{itemize}
Thus by~\cref{lem:lp_bound}, we get that for any sum-intersecting family $\mcF$ of subsets of $[n]$,
$$|\mcF| \le \frac{1}{1 + c_0} 2^n = \frac{1}{1 + 17/8} 2^n = 0.32 \cdot 2^n.$$ \qedhere
\end{proof}
\section{Further work}\label{s:further-work}

With more work, our methods can be extended to give slight improvements to the upper bound we state in~\cref{t:main}; as noted in the introduction, we presented here the simplest version of our arguments here for clarity. However, with the methodology described above, larger and larger finite checks and ad-hoc bounds on probabilities are needed as one gets closer to the conjectured upper bound of~\cref{c:upper}. 

In particular, it would be of interest to understand whether by a more complicated (but finitely supported) choice of $\mu$ in~\cref{lem:lp_bound}, one could prove~\cref{c:upper} via a linear programming bound. 
Tight upper bounds on triangle and $K_4$-intersecting families of graphs~\cite{BZ21,EFF12} have been achieved by similarly inspired methods, but other questions (e.g. the maximum size of a 3-AP-intersecting family) are not immediately amenable to such techniques.

Let us outline the type of computation that we believe could resolve \cref{c:upper} in its entirety, with access to arbitrarily large computational power. The function $\alpha \mapsto |S_\alpha|$ can be expressed in terms of $\left\{\ind_{[\frac 13, \frac 23)}(x\alpha)\right\}_{x \in S}$. Bourgain~\cite{BOU97} had the idea to control the distribution of $|S_{\alpha}|$ by considering the Fourier transform of $\ind_{[\frac 13, \frac 23)}$. Indeed, expanding in terms of the Fourier transform shows that $|S_{\alpha}|$ depends only on the set of linear equations the elements of $S$ satisfy. Moreover, to within a small error, point probabilities of $|S_{\alpha}|$ can be computed given the subset of those linear equations that have coefficients bounded by some explicit function of the error paremeter. As a consequence, the set of feasible distributions for $|S_\alpha|$ (up to controllable additive error) can be reduced to a finite number of cases. Thus, in theory, one could verify an LP bound in the style of~\cref{lem:lp_bound} (after a suitable guess for $\mu$) by a a finite computation; in practice such a computation seemed intractable for us to implement here. 

The following pair of questions naturally extend the problems studied thus far:

\begin{qn} 
What is the size of a largest sum-intersecting family of subsets of $[n]$, $\mcF$, if we require that for each $F_1, F_2 \in \mcF$, the intersection $F_1 \cap F_2$ contains three \emph{distinct} elements $x, y, z \in F_1 \cap F_2$ with $x + y = z$? Are there such ``distinct-sum-intersecting'' families of size larger than $\frac18 2^{n}$?
\end{qn}

Related to a question considered by Bourgain~\cite{BOU97}, the following is also a natural extension of the problem studied here.
\begin{qn}
What is the maximum size of a \emph{3-sum-intersecting} family, i.e. a family $\mcF \subset 2^{[n]}$ such for every pair $F_1, F_2 \in \mcF$, the intersection $F_1 \cap F_2$ must contain $4$ not-necessarily distinct elements $x_1, x_2, x_3, y \in F_1 \cap F_2$ with $x_1 + x_2 + x_3 = y$? Is the family of subsets containing the pair $\{1, 3\}$ (of size $\frac14 2^n$) a maximum size $3$-sum-intersecting family of subsets?

More generally, for any positive integer $k \ge 2$, can we find a \emph{$k$-sum-intersecting family} of subsets of $[n]$ of size $> \frac14 2^n$?
\end{qn}

\bibliographystyle{alpha}
\bibliography{ksat.bib}

\end{document}